\documentclass[11pt]{amsart}
\usepackage{graphicx}
\usepackage{latexsym}
\usepackage{amsfonts}
\usepackage{amssymb}
\usepackage{comment}

\newcommand{\field}[1]{\mathbb{#1}}

\newcommand{\N}{\field{N}}

\newcommand{\R}{\field{R}}
\newcommand{\Z}{\field{Z}}
\newcommand{\myC}{\field{C}}

\title[Flow of complex roots]{MODELING  COMPLEX ROOT MOTION  OF REAL RANDOM POLYNOMIALS UNDER DIFFERENTIATION} \makeatletter

\@addtoreset{equation}{section}
\makeatother
\author{ Andr\'e Galligo }
\address[ Andr\'e Galligo]{Universit\'e de la C\^ote d'Azur, LJAD and INRIA Aromath.\\
Nice, France. }
\email{andre.galligo@univ-cotedazur.fr}


\begin{document}



\renewcommand{\shortauthors}{Galligo}

\begin{abstract}
In this paper, we consider nonlocal, nonlinear partial differential equations to model anisotropic dynamics of complex root sets of random polynomials under differentiation. These equations aim to generalise the recent PDE obtained by Stefan Steinerberger (2019) in the real case, and the PDE obtained by  Sean O'Rourke and Stefan Steinerberger (2020) in the radial case, which amounts to work in 1D. These PDEs approximate dynamics of the complex roots for random polynomials of sufficiently high degree $n$. The unit of the time $t$ corresponds to $n$ differentiations, and the increment $\Delta t$ corresponds to $\frac{1}{n}$. The  general situation in 2D, in particular for complex roots of real polynomials,  was not yet addressed. The purpose of this paper is to present a first attempt in that direction. We assume that the roots are distributed according to a regular distribution with a local homogeneity property (defined in the text), and that this property is maintained under differentiation. This allows us to derive a system of two coupled equations to model the motion. Our system could be interesting for other applications. The paper is illustrated with examples computed with the Maple system.
\end{abstract}


\maketitle

\section{INTRODUCTION}
The analysis of the relation between the zero set of a polynomial  and
the zero set of its derivative goes back to the Gauss-Lucas theorem: the zero set of the derivative lies in the convex hull of the zero
set. The consideration of patterns of roots and  gaps between the roots of a real-valued polynomial under differentiation goes
back to the beginning of the last century. Ten years ago, thanks to Computer Algebra experiments we pointed out in \cite{Gal1}, intriguing patterns created by iterated differentiations of large degrees real random polynomials, and tried an electrostatic interpretation. These phenomena still remain unexplained although they bear resemblance with collective motions, a subject extensively studied, in mathematics, physics and other disciplines, See e.g.the bibliography of \cite{Carrillo}. However important progress has been made in several directions.

\subsection{Pairing} As we experimented in  \cite{Gal1}, the critical points and roots of random polynomial appear to pair with one another and form dotted lines (see examples in section 2).  Hanin investigated the pairing phenomenon between roots and critical points for
several classes of random functions \cite{ha1,ha2,ha3}, including random polynomials with independent roots. He proved that the distance between a fixed, deterministic root and its nearest critical point is $O(\frac{1}{n})$, in the case where $\mu$ has a bounded density supported on the Riemann sphere. Following a conjecture of Pemantle and Rivin, \cite{Pemantle}, Kabluchko and Seidel \cite{Kab}, then O'Rourke and Williams \cite{ORou}, determined the asymptotic fluctuations of the critical point  nearest a given root. They provided a fine stochastic analysis of the local situation and a bound on the Wasserstein distance between the root sets  and critical point sets. So they could combine informations at mesoscopic and at macroscopic scales. 

\vspace{.1cm}

These aspects on the distribution  of roots and critical points are now an active field of research, see the above cited articles and their bibliography, with many nice results and technical precisions. Notably, 
denoting by $X_j$, with$j=1..n$ the roots of the random polynomial $P_n$ such that  the corresponding empirical measure
$$ \mu_n:= \frac {1}{n}\sum_{j=1..n} \delta_{X_j}$$
tends to a measure $\mu$ when $n$ tends to infinity;
denote similarly  $\mu_{n-k}$  the empirical measure attached to the complex roots set of the $k^{th}$ derivative of $P_n$.
It was proved that, in several interesting random contexts, for a fixed integer $k$,  $\mu_{n-k}$ also tends to the same measure $\mu$.
It is expected to be so even if $k(n)=o(n)$  tends to infinity. However, if the differentiation operation is repeated $k=round(tn)$ times, where $t$ is viewed as a ``time" parameter ranging between 0 and 1, some  non-trivial macroscopic dynamics of root sets could show up and help explain the mesoscopic behaviors.

\subsection{All roots are real case}  Several authors  recently investigated the case where all the roots (of the iterated differentiated polynomials) lie on the real axis or on the unit circle.  In 2019, Steinerberger \cite{Stein1} derived the following partial differential equation to describe the asymptotic evolution of  root sets of random polynomials on the real axis, when their degree $n$ tends to infinity:
$$\partial_t u = \frac{1}{\pi}  \partial_x (\arctan (\frac{u}{Hu}))$$
where $u$ is the (regular) limit density of the roots set, $Hu$ is the Hilbert transform of $u$ and $ \partial_x (\arctan (\frac{u}{Hu}))= - \partial_x (\arctan (\frac{Hu}{u}))$. He provided a  rather informal but inspiring construction  of his PDE: following the classical electrical interpretation of a critical point of $P_n$ as an equilibrium of repulsion-attraction forces from the roots of $P_n$, he divided them into a local near field (with a local uniformity property assumption) and an averaged  far field  estimated via a Cauchy-Stieltjes integral, hence a Hilbert transform of the density.
Then, several articles \cite{Gra,Kis, Kab21,Alazar} successively offered more and more complicated and detailed analysis in the periodic setting (i.e. on the circle), they provided rigorous proof of “crystallization” under repeated differentiation. They show global regularity and exponential in time convergence to uniform density. In \cite{Kis} the global in time control follows from the analysis of the propagation of errors equation, with  nonlinear fractional diffusion inspired by similar developments in the modelisation of collective motions.  While in  \cite{Kab21}, the more ''explicit`` strategy of proof relies on complex analysis  consequences of the (proved) fact  that the roots of the iterated  derivatives are distributed according to the free multiplicative convolution of  $\mu$ and a free unitary Poisson distribution. This last work can be related to \cite{tao}, \cite{Stein2},  \cite{HKab21} and \cite{Arizmendi}  for free probabilities and also to \cite{Shapiro} for the saddle point technique of proof.

\subsection{Rotationally invariant case} Assuming that the limit root set is rotationally invariant, the problem reduces to the radial evolution. O'Rourke  and Steinerberger  \cite{ORou-Ste} then Hoskins and Kabluchko \cite{HKab21} proved that a PDE similar but simpler than Steinerberger PDE, could model the asymptotic root motion, reducing the problem to a 1D setting.
In both articles, the previous outlined strategy is again considered but the local near field is reduced to the single nearest root action. For the mean field action, 
O'Rourke  and Steinerberger rely on the following two observations. By symmetry, the velocity of the motion is oriented towards the 
origin. Summing the density on a family of concentric circles of radius $r$, the mean field action is approximated by a sum of integrals parametered by $r$, of the form $\frac{1}{2\pi} \int_0^{2\pi}  \frac{d\theta}{r-s e^{i \theta}}$ which value is either $0$ or $\frac{1}{r}$.
Then, they derive the following PDE:
$$\partial_t \psi =   \partial_r [(\frac{1}{r} \int_0^r \psi)^{-1} \psi ]$$
where $\psi(r,t)=2\pi r u(r,t)$ is the sum of the density $u(r,t)$ on the circle of radius $r$. In that setting the density only depends on the radius. 

Hoskins and Kabluchko \cite{Kab} study the same PDE and provide many interesting quantitative precisions and examples satisfying the radial assumption. They exploit the nice analytic properties of the root sets of the following family of random polynomials, already studied in \cite{K-Z} then by R. Feng and D. Yao \cite{Feng}, 
$$ P_n = \sum_0^n  \xi_k f_{k,n} z^n \;; \;  f_{k,n}= e^{-n v(\frac{k}{n}) +o(n)}
$$
where $\xi_k$ are i.i.d. random variables, e.g. normal real or complex Gaussian ones; and $v:[0,1] \rightarrow \R$ is a continuous function.
Their work and their examples support the correctness of the method.

\subsection{Non rotationally invariant case} The general situation in 2D, for complex roots (without the assumption of radial distribution, which amounts to work in 1D)  is a difficult but important problem since it contains the case of real random polynomials which roots are either real or conjugated complex pairs.    To our best knowledge it was not yet studied (Z. Kabluchko, in the introduction of his very recent article \cite{Kab21}, considered the question but quickly  focused on the 1D case).  The purpose of this paper is to present a first attempt in that direction. 

Up to a translation anda homothety, the sum of the roots is $0$ and they are contained in the unit disk. We use polar coordinates to represent  a complex number or variable $z=\rho e^{i\theta}$.
Our proposition stems on four ideas: 
\begin{enumerate}
 
\item   From our experiments in \cite{Gal1} on the root sets of real random polynomials, we observed that the motion under differentiation is also attracted by the real axis. To capture this turning behavior, we introduce some anisotropy in the model.

\item Steinerberger approach with a local ''repulsive`` field, based on stochastic properties to prevent collision between roots of the same random polynomial, and an ''attractive`` mean field to keep the cohesion of the motion, should be respected. In order to generalise in 2D his local uniformity assumption we, inspired by Lagrangian representations in fluid mechanics, view it as an infinitesimal histogram with balanced bins. 

\item Since the usual notion of histogram in 1D is based on the total ordering of $\R$, we precise our choices for defining an adapted notion of 2D histogram. We borrow from O'Rourke  and Steinerberger the method,  in the radial case, of summing the densities on concentric circles.  Since, in our setting, these densities are not constant on the circles,  we consider a  ''weighted radial  marginal`` distribution $\psi(\rho)$ in polar coordinates of the density $u$;  $\psi(\rho,t)= \int_0^{2\pi}  \rho u(\rho, \theta)$.
We can repeat the argument used in the 1D case, and assume that locally this marginal distribution is constant hence is approximated by evenly spaced points (i.e. the roots lie on evenly spaced circles). We denote by  $a(\rho)$ the spacing.

Now, on each of these circles of radius $\rho$, the roots follow a conditional distribution which  is again assumed locally constant hence is approximated by evenly spaced points, near the point of polar coordinates 
$(\rho,\theta)$.   We denote by
$b(\rho, \theta) =c(\rho, \theta)a(\rho)$ this spacing. Finally, replacing near each root this ''infinitesimal`` portion of ring by a rectangle, we generalise the assumption of Steinerberger in 1D and consider that the near field of a root  is represented by the action of near roots organised in  a $(a, b)$-bi-periodic lattice. Equipped with this geometric formalism,  we can adapt the methodology followed in \cite{Stein1}. 

\item We replace the use in \cite{Stein1}  of the two real trigonometric  functions $cotan(x)$ and $\arctan(x)$ by two adapted complex
functions $F_c (z)$  and $G_c (z)$, depending on the parameter  $c$, such that 

$$G_c (\frac{1}{F_c (z)})=z  \; and \; G_c(z)=z+o(z).$$

Intuitively, the variations of $c$ aim to capture the anisotropy of the trajectories. 
With the Cauchy transform $S_u$, which replaces in 2D the Hilbert transform,  the introduction of $\psi(\rho)$, and with the relations 
$c= \frac{\psi^2}{u} $ (that we establish),  we can compute $G_c (z)$. Then, we will be able to
 derive a system of PDEs  modeling the motion. 
\end{enumerate}
The resulting model consists of the two following relations for the velocity of a root $\xi$ :  
$$v=  \partial_t \xi =  \frac{1}{u} \frac{1}{ b e^{- i \theta}} G_c( b e^{- i \theta} \frac{u}{S_u}))$$
$$  \partial_t b(\rho,\theta,t) =  \frac{b}{\psi}  \partial_t  \psi -  \frac{b}{u}   \partial_t u .$$

Then, the density is modeled by a  PDE, involving a divergence operator:
$$ \partial_t u(\rho,\theta,t) =\nabla . ( \frac{1}{b e^{- i \theta} } G_c(b e^{- i \theta}  \frac{u}{S_u})) \; .$$

\vspace{.2cm}

The paper is organised as follows. Section 2 is devoted to generalities recalled for the convenience of the reader and presents examples. Section 3 is devoted  to  the construction of our 2D histograms  with small rectangles, useful to express our intuition, and to motivate the introduction in  Section 4 of our  main assumption on a local homogeneity property of the distributions. Section 4 also gives insights on the motion of our infinitesimal rectangle and expresses our Lagrangian viewpoint. Section 5 defines the bi-periodic functions  $F_c (z)$  and $G_c (z)$, substitutes of the usual  $cotan(x)$ and $\arctan(x)$ and computes Taylor approximations of them. Section 6 recalls the basic relation between complex roots and critical points of a polynomial and proceeds to the generalisation of Steinerberger decomposition strategy in the context that we have carefully prepared. Section 7 gives the equations of motion  of our model for the dynamics of complex root sets of random polynomials under differentiation. The conclusion discusses limitations of our model and suggests some extensions for future works.

\section{Generalities and examples}
\subsection{Generalities}
 1.  Up to a translation, the sum of the roots of a polynomial $P_n(z)$ can be assumed $0$;  this feature is equivalent to the vanishing of the second coefficient of $P_n$, hence is stable when we replace $P_n(z)$ by all its derivatives. 
 Similarly, if for $\lambda \neq 0$, we consider $Q_n(z)=P_n(\lambda z)$, then $\frac{d}{dz}Q_n(z)= \lambda \frac{d}{dz}P_n(\lambda z)$; so the same homothety-rotation transforms all the root sets of all derivatives of $P_n$. Therefore we can restrict ourselves to the case where the root set  is included in the unit disk.  As a consequence we expect  a scaling of the spacings between the complex roots of order
 $O(\frac{1}{\sqrt{n}})$. Note that if there is an anisotropy in the distribution of the root sets, then it remains after a rotation.
 
 2.  Gauss-Lucas theorem is a consequence of the fact that the logarithmic derivative of $P_n$ vanishes at a critical point $\eta$, when all roots 
 $\{X_j\}$ of   $P_n$ are distinct, which is always assumed for random polynomials, $ \sum_j \frac{1}{\eta - X_j } = \; 0.$ Let $X_1$ be the root nearest to  $\eta$, then
 $$    \frac{1}{\eta - X_1 } =- \sum_{j \neq 1}  \frac{1}{\eta -X_j} \approx - n S_u (\eta) := - n \int_\myC  \frac{u(\zeta) }{ z- \zeta} d \zeta 
 $$
where $S_u$ is the Cauchy transform of the density function $u$. The intuition behind the pairings, quoted in the introduction, departs from this formula, since it shows that (if $S_u(\eta)$ does not vanish)  $| \eta - X_1 | $ is at least smaller than $O(\frac{1}{n})$, hence much smaller than  the  expected order of   distances between the roots, presumed $O(\frac{1}{\sqrt{n}})$.

3. By  Gauss-Lucas theorem, we expect that, under iterated differentiation,  the root sets shrink toward the origin. When the coefficients of  $P_n$ are real, we also expect that, under iterated differentiation, the pairs of conjugated complex roots are also attracted by the real axis.
 This is observed in the computed example at mesoscopic scale. 
 
\subsection{Examples of trajectories}
 We first consider the case of the real root sets of the iterated derivatives of a (real) polynomial $f(x)$ of degree $n$. As we did in \cite{Gal1} we store  the real root sets of the $f^{(k)}$   in a 2D diagram,

 $  \cup_k \; \;\{fsolve(f^{(k)},x)\} \, \times \, \{n-k\} .$
 
The picture at the  top of Figure 1 corresponds to a polynomial  $f(x)$  with an all real root set. The bottom pictures corresponds to a polynomial  $f(x)$ which has also complex roots; under differentiation they ''fall`` on the real axis and appear as new ''curl`` of trajectories in the diagram.
 In the right picture, we relied on  fractional derivatives, which interpolate the obtained dotted curves of the trajectories, to clarify the motion indicated by decreasing degrees.

\begin{figure}[ht!]
  \includegraphics[width=14cm]{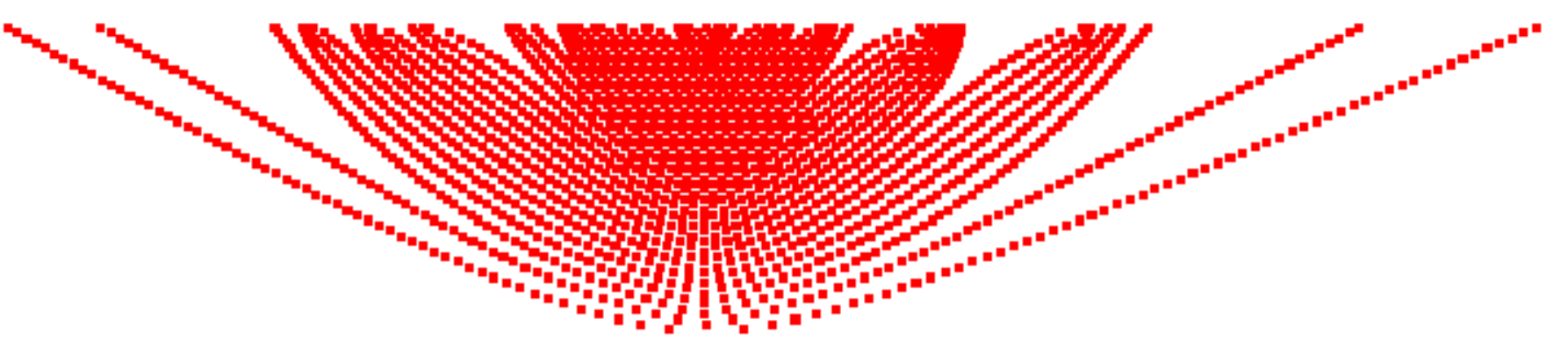}
  \includegraphics[width=4cm]{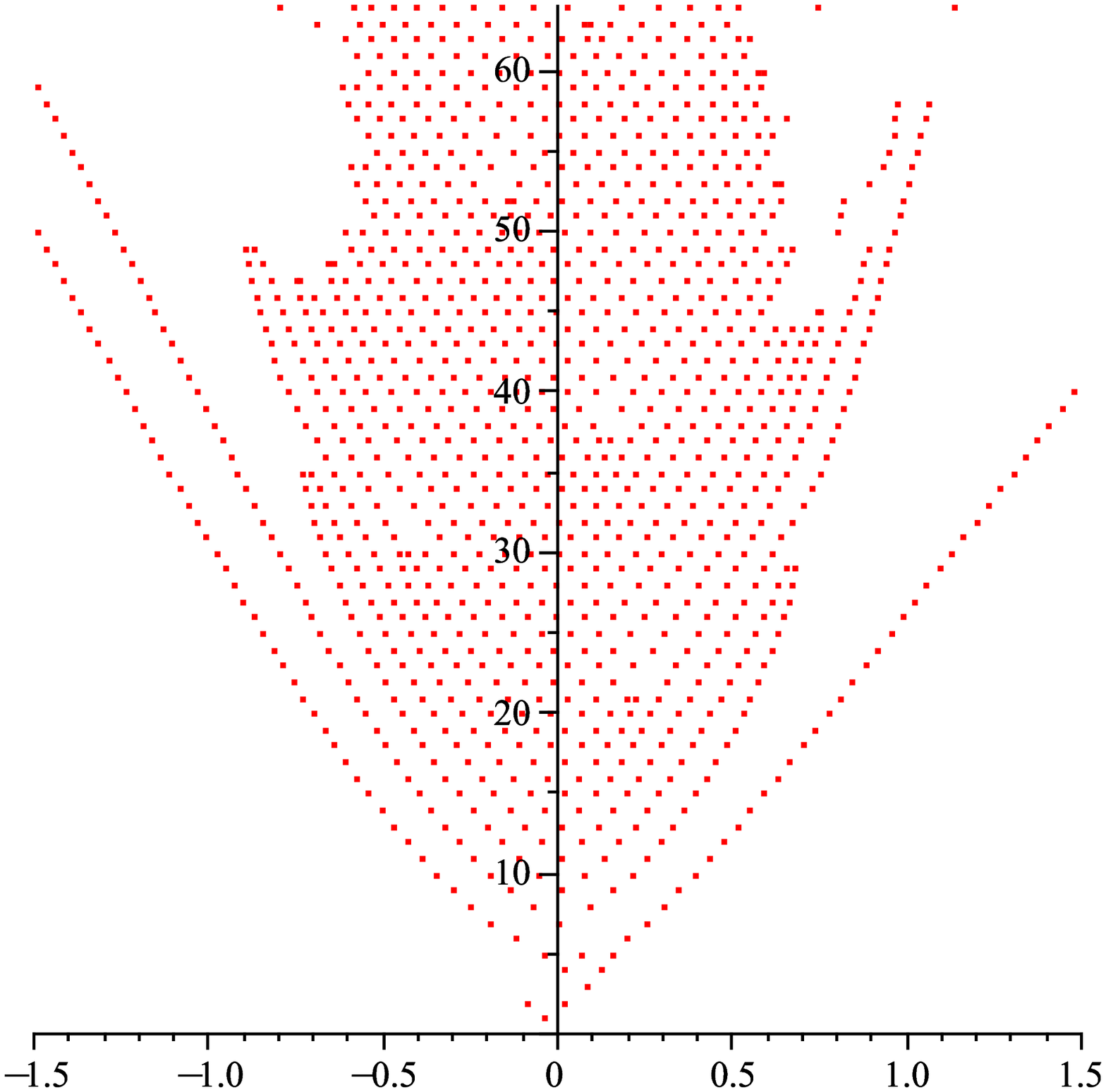}
      \includegraphics[width=4cm]{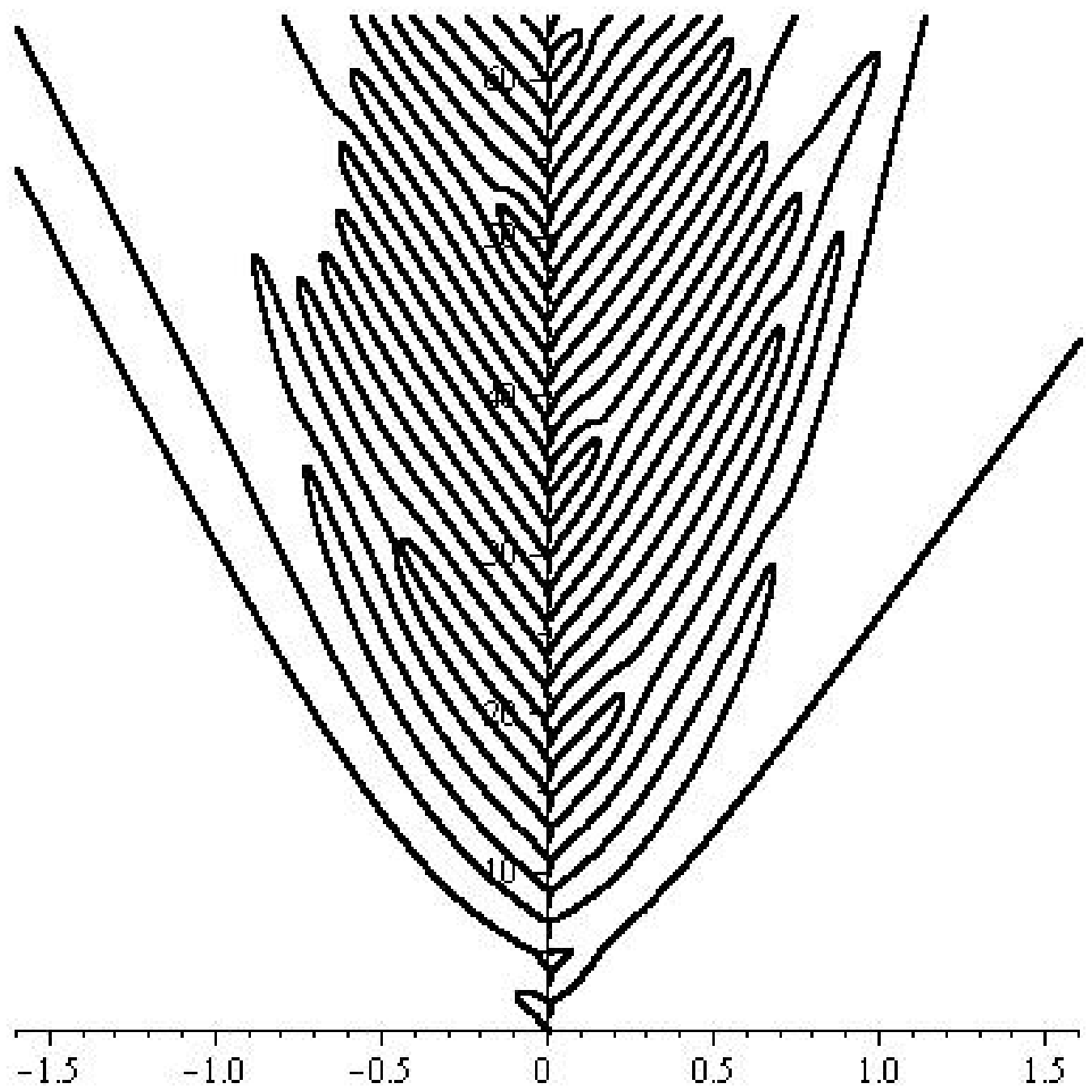}
 \caption{In the real setting}
\end{figure}  

2. For complex root sets we represent the decreasing degrees, showing the motion, by colors:  in Figures 2 and 3 the initial 150 random roots are indicated by black solid boxes, then we performed two steps of 35 derivations and colored the intermediate root sets in red, then blue, then pink, then green. The proportion of initial real roots in Figure 2  is about a tenth. The proportion of initial real roots in Figure 3 is about a third, so the attraction towards the real axis is stronger and one can figure out the turning trajectories.

\begin{figure}[ht!]
 
      \includegraphics[width=8cm]{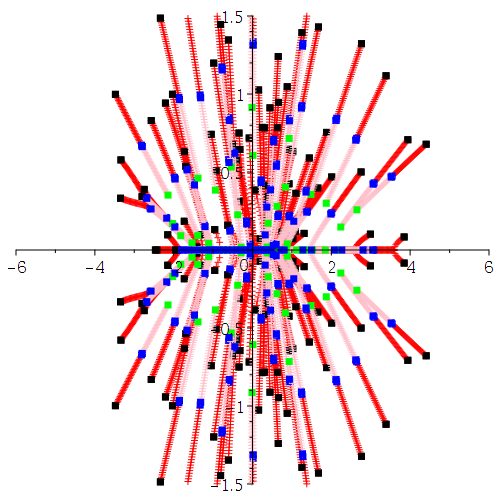}
     \caption{In the complex setting}
\end{figure}  

\begin{figure}[ht!]
 
      \includegraphics[width=8cm]{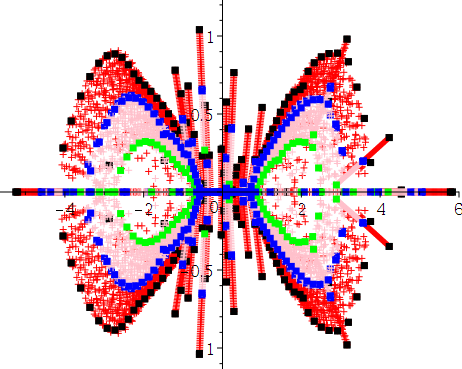}
     \caption{Turning trajectories}
\end{figure}

\section{2D histograms}
In this section, we present a mesoscopic construction which guides our intuition for describing, in the next section, our main hypothesis and our choice of parameters at the macroscopic scale. 
\subsection{Homogeneity hypothesis}
When all the $n$ roots of a random polynomial $P_n(x)$, with $n\gg1$, are real, the empirical distribution of its roots can be approximately represented by a histogram with $m$  bins $I_k$
such that
$Prob(I_k) =\frac{1}{m}$. For a positive integer $s$  and for large $m$,  with $n=ms$, the root set of  $P_n(x)$ can be approximated by the collection made by evenly spaced $s$ points in each of the $m$ bins. 
Then the local interaction between the near roots is reduced to the case of a lattice. 
We view this construction as one of the key observation in the derivation of Steinerberger PDE in \cite{Stein1}, and we would like to generalise it for complex roots sets.

However, this simple construction, of  moving bins, relies on the fact that $\R$ is naturally totally ordered, which is not true for $\R^2$.
\subsection{Concentric rings}
A random complex variable $z$ represented by
$z= \rho e^{i \theta}$ follows  the joint law of the two real random variables $\rho$ and $ \theta$. To mimic the construction of the previous subsection, we will first construct $m^2$ rectangular bins forming a 2D histogram of the uniform probability distribution on the unit disk, where we will distribute ''uniformly`` $n=m^2 s^2$ points, with $m\gg1$ and $s\gg1$ . These  rectangles approximate tiny sectors of concentric rings with an area equal to $ \frac{\pi}{m^2}$. For a general probability distribution $u(\rho, \theta)$, with support in the unit disk, with zones where the distribution vanishes, there will still be  $m^2$ rectangular bins, approximating tiny sectors of concentric rings; but their locations and areas will depend on  $u(\rho, \theta)$ as we now describe.

We define the weighted marginal probability density of $u(\rho, \theta)$  to be  $\psi(\rho)= \int_0^{2\pi} u(\rho,\theta) \rho  d\theta$. Note that 
  $\rho u(\rho, \theta)=\psi(\rho) U_{\rho}(\theta) $ with $\int_0^{1} \psi(\rho) d\rho  = 1$ and $ \int_0^{2 \pi} U_{\rho}(\theta) d \theta = 1$. Then, following   \cite{Stein1}, we assume that these 1D probability densities are locally well approximated, on their supports,  by constant values. Consequently the roots  will be suitably equi-distributed into each of the $m^2$ bins.
  
  In other words, for large integers $s$  and  $m$,  with $n=m^2s^2$, we consider $m$ disjoint rings $T_k$, with  $k=1..m$, each bordered by two concentric circles, and containing  $m s^2$  roots for a total of $n$ roots. There might be void rings between the $T_k$.
Then, we subdivide the support of $u(\rho, \theta)$ in each ring $T_k$ into $m$ sectors 

$S_{jk}$, $j=1..2k-1, k=1..m$, with  
$Prob(S_{jk}) =\frac{1}{m^2}$. 

When $m$ is very large, we approximate the  $S_{jk}$, by small rectangles and assume that they contain $s^2$ roots   organised in a lattice  evenly spaced in the two orthogonal directions $d\rho$ and $d\theta$. We denote by $a(\rho)$ the spacing in the $d\rho$ direction, and by  $b(\rho, \theta)=c(\rho, \theta). a(\rho)$ the spacing in the $d\theta$ direction. We have the  following estimations
$$ a(\rho) = \frac{1}{\sqrt{n} \psi(\rho)} \;; \;  b(\rho, \theta) = \frac{\psi(\rho)}{\sqrt{n} u(\rho, \theta)} \;; \;  c(\rho, \theta) = \frac{\psi(\rho)^2}{u(\rho, \theta)}.
$$
Then, as expected $a b u = \frac{1}{n}$, but also  $ \frac{1}{n}=a^2 \psi^2 = \frac{b^2 u^2}{\psi^2}$.

Figure 4 sketches,  approximately  and at mesoscopic scale, a rectangle (in blue), with $s^2$ bi-periodically  spaced roots (in black); and the transformed rectangle (in red) after performing a differentiation. For visibility, the effect has been strongly  accentuated.  
\begin{figure}[ht!]
      \includegraphics[width=9cm]{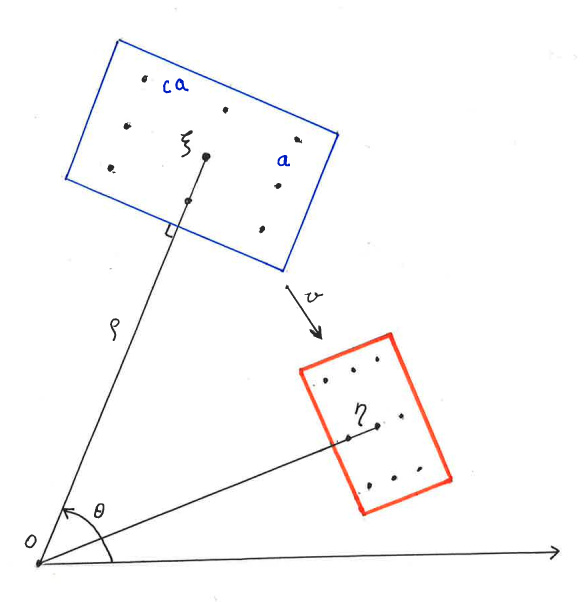}
   \caption{Moving rectangles}
\end{figure}  

\section{Main assumption and useful formulas}
The previous construction motivates the following assumption on the density functions $u(\rho, \theta)$ of the limit empirical measures of  the roots of polynomials $P_n$,   when $n$ tends to infinity. This property is also assumed to hold after iterated differentiations.

\subsection{Homogeneity property for $u(\rho, \theta)$}
The density probability $u(\rho, \theta)$  is smooth; its support is included in the unit disk of $\myC$. 
Near any point of its support, the density is locally the product of the two probability densities $\psi(\rho)$ and  $U_{\rho}(\theta)$, which are also assumed smooth; 
$\psi(\rho)= \int_0^{2\pi} u(\rho,\theta) \rho  d\theta$. 
Moreover, with  $c= \frac{\psi^2}{u}$, the local density  $u(\rho, \theta)$ can be approximated by the empiric density of a lattice  
 $ \rho e^{i \theta}+ a e^{i \theta}(\Z + i c \Z)$  included in  $\myC$.

\subsection{Motion of rectangles}
In this subsection, we consider the motion of the root set, organised in small rectangles as described in the previous section, with respect to an abstract time $t$. 

The rectangles are homothetic to a rectangle 
 $R(\xi,a,ca,\theta,t)$, which is characterised by its center $\xi=\rho e^{i \theta}$, its direction $\theta$ and a
 corner $\xi_1= \xi + a e^{i \theta}(1+ i c)$. 
 The infinitesimal motion of the rectangles, corresponds to a composition of classical transformations:  translation, rotation, dilatation, elongation. It is described by the time derivatives of the functions $(\rho(t), \theta(t), a(t), b(t), c(t))$.
 The difference  $\xi_1 -\xi$ can be viewed as a spatial increment of  $\xi(t)$. Differentiating  $\xi(t)$ and $\xi_1(t)$ with respect to the time $t$  and taking the limit when this spatial increment tends to $0$ will be used below to derive a mass conservation relation,
needed in Section 7 to derive the equations of motion.

{\bf Remark} Heuristically, the transformations of the rectangles $R(\xi,a,ca,\theta,t)$ during the motion capture the turning behavior of the root set flow. The rectangles move towards the origin together with the ring $T_k(t)$ which contains them;  and they can also move towards the real axis depending on the distribution $c(t)$.

It would be worthwhile to compare the capabilities of description of a flow by moving infinitesimal rectangles  with other formalisms, notably, the recent elaborate Hamiltonian model developed by A. Cavagna et al. in  \cite{Ham} to study the rotations of bird flocks, or the recent tensor formalism developed by B. During et al \cite{During}, to model finger prints creation.

 In the last section, we will discuss weaker versions of our homogeneity property where the rectangles centered at $\xi$ are not aligned with $\xi$. The aim is  to  precisely capture the directions of ''stress' created by the collective motion.

\subsection{Conservation of mass and material derivatives}

In an Eulerian representation, the mass conservation equation for a fluid expresses that
the flux going into a control volume is subtracted from
the flux going out of a control volume. 
The limit is taken as the control volume shrinks to a point, and is computed with a
divergence operator. 

However, the Lagrangian formalism, which seems better suited for our study, presumes to know the starting position of a particle
and treat the subsequent position as a dependent variable
(while in the Eulerian representation, this is treated as the independent variable).
Another time derivation, often called material derivation, is denoted by $\frac{D}{Dt}$ or simply by $D_t$ such that 
for a field $f$, denoting by $ \nabla f $  the gradient of $f$ and by $ \vec{v}$  the velocity vector,
$$D_t f := \partial_t f +  \vec{v} .  \nabla f .$$

Then the mass conservation is written in two different ways, having the same meaning:
$$  D_t u = -u \nabla . \vec{v}; \;  \;  or  \; 
  \partial_t u = - \nabla . u \vec{v} = -u \nabla . \vec{v} - \vec{v} .  \nabla u.$$

\subsection{Useful formulas}

In our 1D setting, we have the relation $ u a = \frac{1}{n}$, which implies $  D_t u = - \frac{u}{a} D_t a$. 
Applying the equilibrium formula, observed by Steinerberger,
$$ \pi n u . cotan(\pi n u (\eta-\xi)) = - \pi n H_u$$
solving it in $\eta-\xi $ then dividing by $\Delta t = \frac{1}{n}$, 
one obtains $$ v= - \arctan(\frac{u}{Hu})\frac{1}{\pi u}.$$

The material time derivative of the spacing  $a$ is obtained by applying the previous relation  to two near-by points $\xi$ and $\xi_1$,  and subtracting, we get with obvious notations
$$
 \frac{(\eta_1 - \eta) - (\xi_1 - \xi)}{\Delta t}=  \frac{(\eta_1 - \xi_1))}{\Delta t}-  \frac{(\eta - \xi))}{\Delta t}
$$
such that the left hand side approximates $D_t a$; while the the right hand, side once divided by the spacing $a$, approximates $ \partial_x v$. Hence  $$D_t a = a \partial_x v .$$

This implies  $D_t u = - u \partial_x v$  and we recover Steinerberger PDE, 
 $$ \partial_t u = \frac{1}{\pi} \partial_x (\arctan(\frac{u}{Hu})).$$

\noindent
In 2D, we proceed similarly but with the two spacings $a$ and $b$ in two orthogonal directions. We now have 
$$ D_t u = -( \frac{u}{a} D_t a + \frac{u}{b} D_t b);
$$
We will compute the velocity in Section 6. It will be a vector $\vec{v}$ represented as a complex number.
Repeating twice the reasoning  we just made, with a pair of near-by roots spaced by $a e^{i \theta}$, respectively by $i b e^{i \theta}$,
we get the two components of the velocity (in the orthogonal frame $(\xi, X,Y)) $ in the direction $\theta$) denoted by $v_X$ and $v_Y$; then
 $D_t a = a \partial_X v_X$ and $D_t b = b \partial_Y v_Y$. Hence,  
$$u  \nabla . \vec{v}= u (\partial_X v_X + \partial_Y v_Y )=  \frac{u}{a} D_t a +  \frac{u}{b} D_t b = -D_t u .$$
Translating in Eulerian coordinates, we finally get the familiar equation,
$$  \partial_t u = - \nabla . u \vec{v}.
$$

\section{Bi-periodic functions}

For $c>0$, we first consider the lattice $\Lambda(c) = \Z + i c \Z$ included in $\myC$, together with the odd bi-periodic function $F_c (z)$, (we will use the notation $\Lambda(c)^* = \Lambda(c) - \{0\}$);

$$ F_c (z) = \sum_{\lambda \in \Lambda(c)} \frac{1}{z- \lambda} \; = \frac{1}{z} + \sum_{\lambda \in \Lambda(c)^*} \frac{z}{z^2- \lambda^2}. $$

It would be nice to express $F_c$ in terms of classical zeta functions, but it is also useful to  approximate it, near $z=0$, by $\frac{1}{z}$ times a Taylor expansion

$$ \frac{1}{z} (1- z^2 \sum_{\lambda \in \Lambda(c)^*}  \frac{1}{\lambda^2} +  \sum_ {m=2}^{\infty} (-1)^m z^{2m} \sum_ {\lambda \in \Lambda(c)^*}  \frac{1}{\lambda^{2m}}) $$
that we may truncate at order $o(z^5)$ for a computable approximation.

The series $\sum_{\lambda \in \Lambda(c)^*}  \frac{1}{\lambda^{2m}}$ converges absolutely for $m>1$ and is the value of a classical special function.
In the first sum, we group opposite and conjugated terms to get a convergent series, that we denote by $g(c)$, with  $g(1)=0$ and $g(\frac{1}{c} )= -c^2 g(c)$.
For $m=2$, we call it $h(c)$ and note that  $h(\frac{1}{c} )= c^4 h(c)$ and that $h$ does not vanish for $c=1$.
The sums   
$$ \sum_{l \in \N^*}  \frac{1}{l^2} = \frac{\pi^2}{6} \;; \; \sum_{l \in \N^*}  \frac{1}{l^4} = \frac{\pi^4}{90}$$
are well known  (Bernoulli numbers). By a simple computation, 

$g(c) =  \frac{\pi^2}{3}(1-\frac{1}{c^2}) +2 g_1(c)$ and 
$h(c) =  \frac{\pi^4}{45}(1+\frac{1}{c^4}) +2 h_1(c)$, with

$$ g_1(c):=\sum_{(l,k) \in \N^* \times \N^*}  \frac{l^2-k^2 c^2}{(l^2+k^2 c^2)^2};$$ 
$$h_1(c):=  \sum_{(l,k) \in \N^* \times \N^*}  \frac{l^4-6 k^2 l^2 c^2+c^4 k^4}{(l^2+k^2 c^2)^4}; $$
Approximate graph of $g_1(c)$ and  $h_1(c)$,  for $c=0.1..10$, resp. for $c=0.1..3$,  are shown in Figure 5.

\begin{figure}[ht!]
 
      \includegraphics[width=5cm]{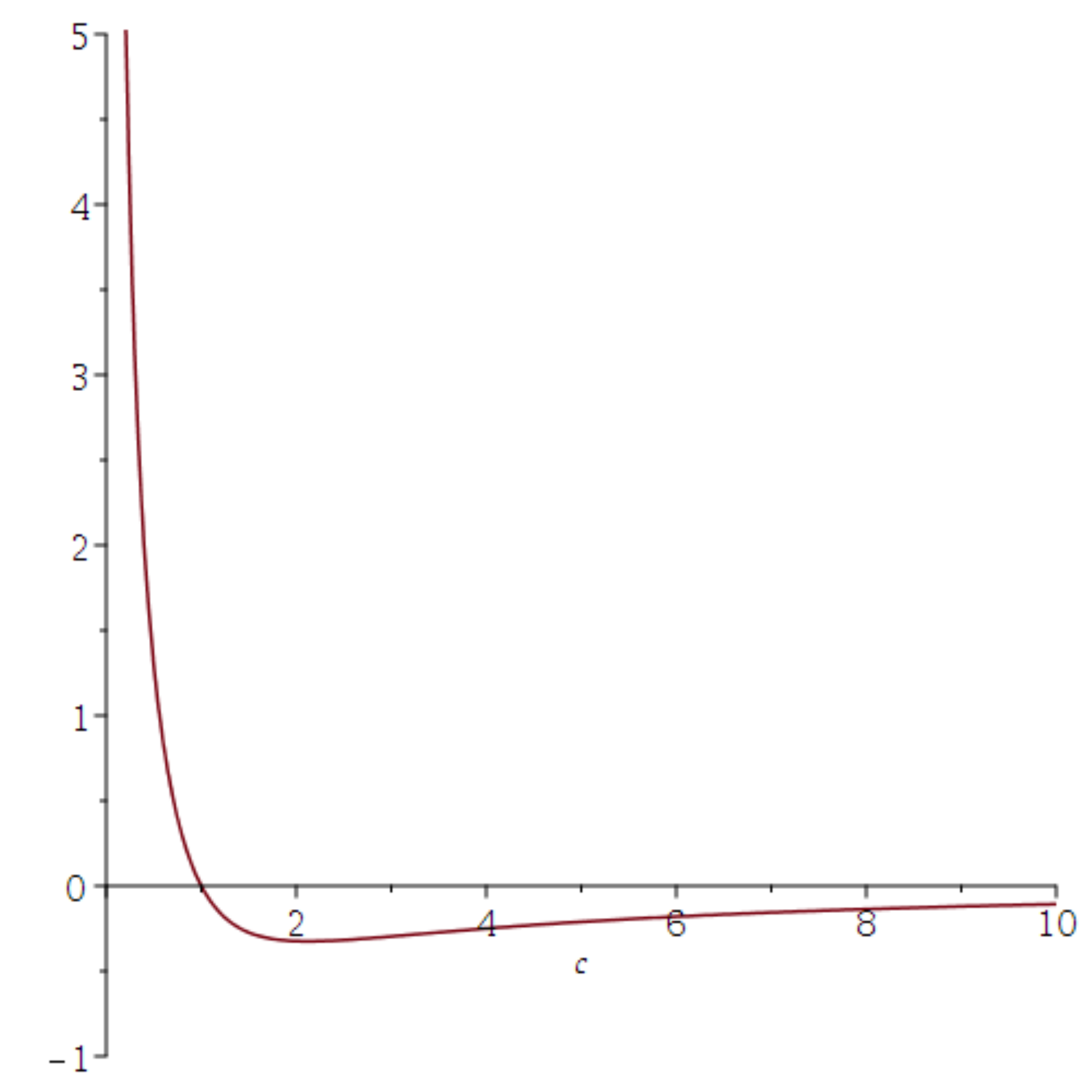}
      \includegraphics[width=5cm]{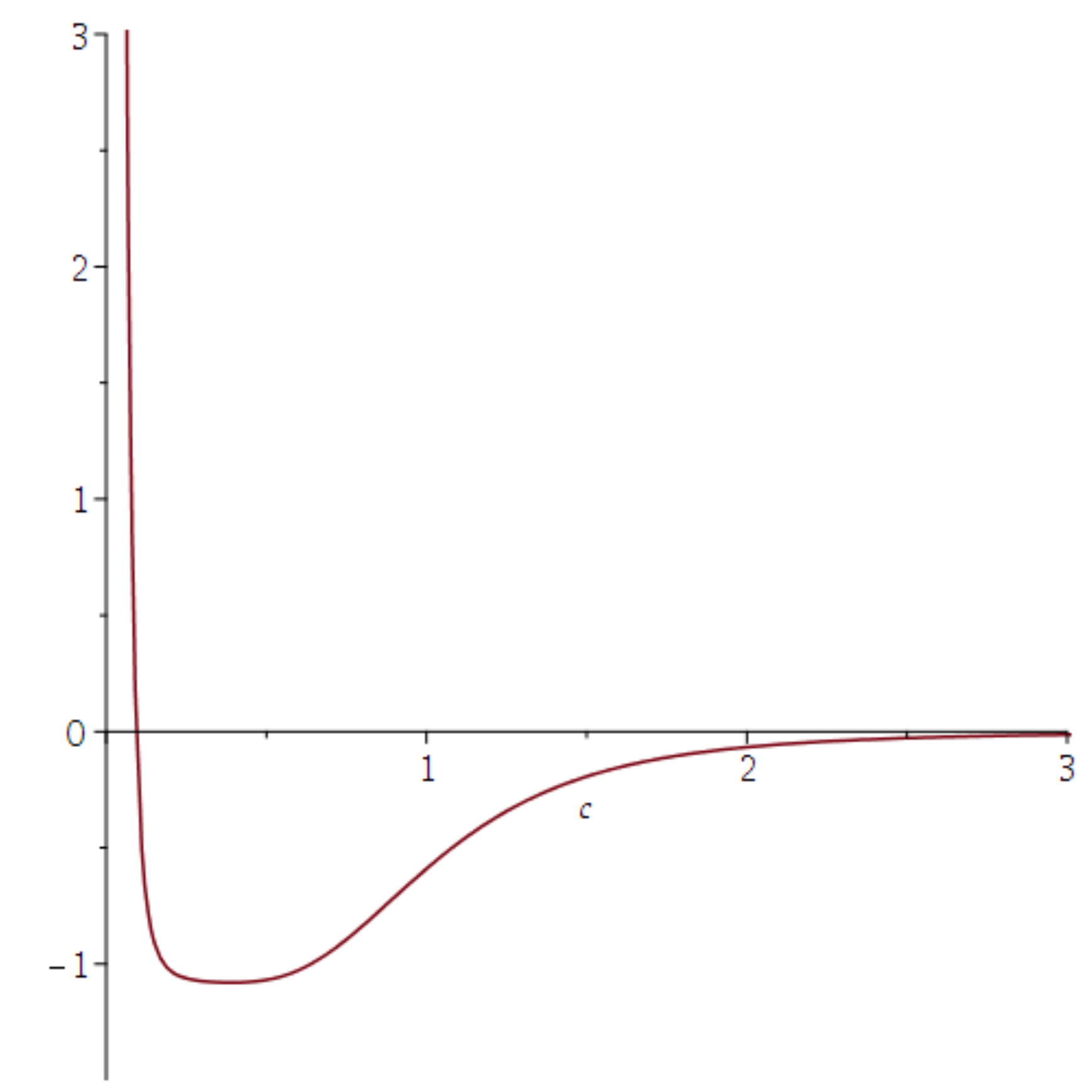}
     \caption{graphs of $g_1$ and $h_1$}
\end{figure}  

\subsection{The inverse function $G_c$}

For each fixed $c>0$, $F_c(z)$ is a substitute for the role played by $\pi cotan(\pi x)$ in \cite{Stein1}. 
We now mimic the  construction of the real function $\arctan(w)$, from the real function $cotan(x)$ by solving the equation  $cotan(x)= \frac{1}{w}$ to get $x=\arctan(w)$.

We denote by $z=G_c(w)$ the solution $z$ of the equation $F_c(z)=\frac{1}{w}$. This defines another complex function
$G_c$ depending on the parameter $c$. To compute a Taylor approximation at order 6 of  $z=G_c(w)$ near $w=0$, we proceed to 
the  inversion on the previous  approximations of  $zF_c(z)$  by a Taylor expansion. We get

$ w = z( 1 - g(c) z^2 + h(c) z^4 + o(z^5) )^{-1} ,$   hence

$ w = z - g(c) z^3 + ( g(c)^2 - h(c)) z^5 + o(z^6) , $
which provides
$$ z = w + g(c) w^3 + (g (c)^2 - h(c)) w^5 + o(w^6) . $$
So approximations of $G_c(z)$ are easily computable.
\subsection{Change of coordinates}

We now introduce the lattice $\Lambda(a, ca, \theta) = a e^{i \theta}(\Z + i c \Z)$ included in $\myC$, corresponding to the rectangle considered in the previous section, together with the bi-periodic function $F_{a,ca,\theta}$ :

$$ F_{a,c,\theta}  (z-\xi) = \sum_{\lambda \in \Lambda({a,ac,\theta})} \frac{1}{z-\xi- \lambda} ,$$

a simple calculation shows that:

$$ F_{a,ca,\theta} (z-\xi) = \frac{1}{ae^{i \theta}} F_c(\frac{z-\xi}{ae^{i \theta}}).$$

Again, we can solve $F_{a,ca,\theta} (z-\xi)=\frac{1}{w}$  and we now get:

$$ z-\xi = a e^{i \theta} G_c(\frac{w}{ a e^{i \theta}}) 
$$

and obviously near $w=0$, $ z-\xi = w + o(w^2)$.

\section{Electrostatic interpretation}

Writing the logarithmic derivative of $P(z)$, we get

$$\frac{P'(z)}{P(z)}= \sum_{j=1..n}  \frac{1}{z-X_j}$$
where the $X_j$ are the $n$ complex roots of $P(z)$. We choose one of these roots, that we denote by $\xi$ and we denote by $\eta$ the root of the derivative of $P(z)$ nearest to $\xi$.
The  asymptotic (when $n$ tends to infinity) relative location and properties of such pairing $\xi, \eta$ has been studied in detail and for a number of random situations by several authors, as related in our introduction; assuming that the empirical measure $ \mu_n := \frac{1}{n}\sum_{j=1..n} \delta_{X_j}$ converges to a measure with a sufficiently regular density $u(z)$ on $\myC$. A key tool for these studies is the Cauchy  transform (also called Cauchy Stieltjes transform) of $u$:
$$ S_u (z) := \int_\myC  \frac{u(\zeta) }{ z- \zeta} d \zeta $$
it is well defined when the integral converges.
Let us give an illustrative example (from \cite{ORou}) of this transform in a simple case.

Assume that the distribution $u$ is uniform in a disc of radius $R$ centered at the origin and take $z=\rho e^{i \phi}$, we have

$ S_u(z)=\frac{1}{z} $ if $|z| \geq R \; $, $ \; S_u(z)=\frac{\bar{z}}{R^2} $ if $|z| < R.$

As noted in \cite{ORou}, heuristically, when $ S_u (z)$ is zero (or almost zero), there is no root of $P^{'}$ near $z$; whereas when $ S_u (\xi)$ is clearly non zero, there is a root $\eta$ of $P^{'}$, such that $|\xi-\eta| = O(\frac {1}{n})$.
Indeed, since $P^{'}(\eta)=0$ we then have:
$$ \frac {1}{\eta - \xi} = - \sum_{X_j \neq \xi}  \frac{1}{\eta-X_j}$$
and the right-hand side is approximated by $-(n-1) S_u(\eta)$ or asymptotically equivalently by $-n S_u(\xi)$.

Steinerberger \cite{Stein1} used the Hilbert transform which, in the case where all roots are real, plays on the real axis the same role as the Cauchy transform on $\myC$.

When $n$ tends to infinity, the previous right-hand side continues to be well approximated by $-n S_u(\xi)$, even if we move $s$ terms, $s \gg n$,  to the left-hand side and write
$$ \frac {1}{\eta - \xi} + \sum_{0<|X_j-\xi|  small}  \frac{1}{(\eta-\xi)-(X_j-\xi))}= - \sum_{|X_j-\xi| large}  \frac{1}{\eta-X_j}.$$

Now we can apply our local homogeneity hypothesis on the limit measure and its expression computed in the previous section, then combine it with  the  equality $n= \frac{1}{a^2 \psi^2}$.
Equating the two approximations, we get the following asymptotic equation:

$$ \frac{1}{ae^{i \theta}} F_c(\frac{\eta-\xi}{ae^{i \theta}}) = -\frac{1}{a^2 \psi^2}S_u$$
where the functions $a, c, \phi, u, S_u$ are evaluated at $\xi$.

This equation can be solved in $\eta-\xi$ near $0$ to get, since $F_c$ is odd,
$$ \eta-\xi = - ae^{i \theta} G_c( \frac{a \psi^2 e^{- i \theta}}{S_u}). $$

To estimate the time derivatives, we need to approximate the quotients by 
$\Delta t :=  \frac{1}{n} =a^2 \psi^2$. So the velocity $v(z,t)$  of a root $z=\xi$ is the limit of $ \frac{\eta-\xi}{\Delta t}$.  
We obtain:
$$v(z,t) = -  \frac{1}{a\psi^2 e^{- i \theta}} G_c( \frac{a\psi^2 e^{- i \theta}}{S_u}) $$

or equivalently, since $bu = a\psi^2 $, the nicer expression:
$$v(z,t) = -  \frac{1}{u}\frac{1}{ b e^{- i \theta}} G_c(b e^{- i \theta} \frac{u}{S_u}). $$

We recalls the relations $c= \frac{\psi^2}{u}$ and   $b = \frac{\psi}{\sqrt{n} u}$. 


With these relations, we are ready to derive the equations of motion of the flow of root sets under differentiation.

\section{Equations of motion}

Here the complex numbers are identified with vectors of $\R^2$ in order to apply the divergence operation  $\nabla . \vec{v}$.  

We can replace in the formula of Subsection 4.4 the expression of $v(z,t)$ established in Section 6.  We obtain the following PDE (with the notation $z=\rho e^{i\theta}$):

$$ \partial_t u(z,t) = \nabla . \left(  \frac{1}{be^{- i \theta}} G_c(b e^{- i \theta} \frac{u}{S_u}) \right).$$

It will be coupled with the following equation obtained  from the logarithmic derivation of  $b = \frac{\psi}{\sqrt{n} u}$,

 $$  \partial_t b(\rho,\theta,t) =  \frac{b}{\psi}  \partial_t  \psi -  \frac{b}{u}   \partial_t u .$$

{\bf Remark 1} As noticed by Zakhar Kabluchko in \cite{Kab21}, to take into account that at each derivation the degree of the polynomial decreases by 1, the density $u(z,t)$ should be the density of the measure $(1-t) \frac{d \mu}{dz}$, where $\mu$ is the measure at time $t=0$. 

\vspace{.2cm}

{\bf Remark 2}
When we assume that the limit distribution of the roots is rotationally invariant and remains so under differentiation, the model can be simplified.

Indeed, on the one hand  the expression of $S_u$ can be explicitly computed following the observations of    O'Rourke and Steinerberger in \cite{ORou-Ste}. On  the other hand, the quantities  $b$ and $c$ will also 
depend only on $\rho$, nevertheless, in that case, our formulation seems more complicated than that of \cite{ORou-Ste}. We will investigate this point in a future work in relation with ``almost'' rotationally invariant situations. We refere to \cite{Feng}  for an analysis of the case of Kac polynomials.  The root sets of a real Kac polynomial are shown in Figure 6.

\begin{figure}[ht!]
 
      \includegraphics[width=6cm]{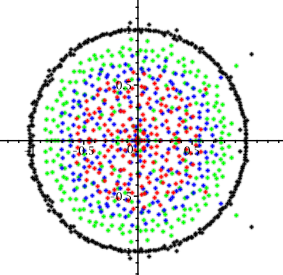}
     \caption{Iterated differentiations of a Kac polynomial.}
\end{figure}  

\section{Conclusion}

In this conference paper, we proposed the first generalisation in 2D of the assumption and the strategy deployed by Steinerberger in 1D to derive its non local PDE in the celebrated article \cite{Stein1}. It  describes the asymptotic evolution of the root set of a random polynomial, under iterated differentiations. We derived a system of two PDEs, the expression of one of the two resembles the  non local non linear Steinerberger PDE. The spatial derivative is, as expected, replaced by a divergence. 

Building from our own interpretation of  the ideas developed in the two articles \cite{Stein1, ORou-Ste}, we constructed an adapted formalism of small sectors of rings (approximated by small rectangles),  moving on evolving concentric rings and which can behave differently in the two orthogonal directions. Our aim was to represent  as  ''simply`` as possible, a collective motion of roots attracted both by the origin and by the real axis. The first attraction is captured by a ''shrinking`` of the concentric rings, while the second attraction is captured by ''sliding`` the distribution inside the ring towards the real axis.  

One such rectangle is characterised by its center $\xi=\rho e{i\theta}$, which represents a root, and by the pair $(\xi+a e{i\theta} , \xi+i b e{i\theta}) $, which represents two near-by roots of $\xi$ in two orthogonal directions. So that $a$ and $b$ can be interpreted as spatial increment in these directions; indeed  these two near-by roots are the generators of a subset of near roots organised in a lattice structure.  

With this restriction,  the rectangle is aligned  with the directions of polar coordinates attached to $\xi=\rho e{i\theta}$, a property designed to simplify the computations.  We notice that during the motion, the velocity of $\xi$ and  $\xi+a e{i\theta}+i b e{i\theta}$ are almost equal. But their tiny difference creates a kind of spin which modifies the features of the moving rectangle. In our formalism, this is captured by the variations of $a$ and $b$, but mostly by that of their quotient $c=\frac{b}{a}$.

We were able to measure it by the formula $c= \frac{\psi^2}{u}$ which plays a key role in our computations. Here $u(\rho, \theta)$ is the limit density function at $\xi$, while $\psi(\rho)$ is a weighted marginal of this density on a circle. This value $c$ also serves as a parameter to construct a function $G_c(z)$,  which generalises in 2D the  trigonometric function  $\arctan(x)$ used in \cite{Stein1} and which can be easily evaluated. 

In 2D, since there are two directions of spatial increment, one needs to consider  two functions $a$ and $b$, instead of only one in 1D. We were able, thanks to the introduction of the marginal distribution $\psi$, to  reduce our study to the case where the spatial increment $a$ does not depends on $\theta$. This feature allowed us to eliminate $a$ from our main equation of motion. However the expression of our main PDE still involves $b(\rho, \theta)$. We could couple this ''mass conservation`` equation  with another PDE, which regulates the evolution of $b$.
\subsection{Discussion and directions for extensions}
- Our formalism allows to describe the asymptotic evolution of some complex root sets which are not invariant by rotations centered at the origin, a limitation imposed by previous related articles. However our main assumption is still restrictive. A weaker assumption  could have been to allow the rectangles to move more freely, by introducing a new variable: the angle  between the vector $\xi$ and one of the sides of the rectangle. This variable would allow to capture more precisely the evolving direction of the local stress created by the collective motion. In an exploratory work, we considered this possibility, but it implied more complicated computations.

- It would be worthwhile to program an efficient solver adapted to our non local PDEs, yet it is not an easy project.

- Then, a first task would be to test this solver with the families of examples studied by O’Rourke and Steinerberger
\cite{ORou-Ste}, for comparing the results.

- A related task would be to proceed to a stability analysis of the radial case by considering limit behaviors of almost radial cases satisfying our homogeneity hypothesis.

- In our model, we assumed a decomposition of the joint density (in polar coordinates) $ u =  \frac{1}{ \rho} \psi( \rho) U( \rho,  \theta)$; therefore the model is limited by the possible explosion at a time $ T  \geq 0$. We expect that  $T$ is positive, but we should prove it rigorously. This explosion may happen if the support of the limit density shrinks to a curve.

- Therefore a suitable setting, able to deal with many more geometric situations, is to assume that the support is made of a domain $D_1$ union a curve $C_1$. To apply the strategy of computation described in this paper, we will assume our homogeneity hypothesis in the interior of  $D_1$, Steinerberger’s hypothesis on the curve and a mixed one on the border of $D_1$. Then, we could estimate the corresponding local contributions as we did in Section 6 and derive the equation of motion as in Section 7.

- So, our present work can be seen as an important first step towards this more general project. We plan to develop this project, starting from the case when the curve is known in advance e.g. if it is the real line or the unit circle.

- One of the anonymous reviewers wrote that the alignments of the root configurations reminded him the works of Boris Shapiro and his collaborators  \cite{Shapiro1, Shapiro2, Shapiro3} on classification of eigenpolynomials for an exactly solvable rational differential operator, and the study of their roots. This similarity also appears with another article of the same team \cite{Shapiro} used by Kabluchko in his recent paper \cite{Kab21} on root sets of iterated derivatives of trigonometric polynomials $f( \theta)$ on the unit circle. In that paper, the derivation of $f( \theta)$ w.r.t. $ \theta$ is interpreted as the action of the differential operator 
$ i(z \frac{d}{dz} - \frac{n}{2})$ on a complex polynomial.

- This suggests that our framework could be extended to modeling polynomial complex root sets under iterative action of specific differential operators. We will develop this point in a future work.

- It would be great to obtain, with our kind of hypothesis, a Hamiltonian formulation which would take advantage of the symmetry of the interactions between the roots at any given time. An inspiring work in that direction is \cite{Ham}, where the authors provide an original Hamiltonian model for the  natural  flocking of birds with adapted internal variables.

- Another direction of research could be to concentrate on special families of random polynomials where one could express and exploit a potential for collective interaction between the roots, similar to the one existing for characteristic polynomials of matrices in the GUE ensemble.  

- In the same spirit, it would be interesting to compare our formalism with  recent developments studying compression of Coulomb gas in 2D or 3D, see e.g. \cite{Butez,Chafai} and the excellent presentations of Sylvia Serfaty at the last mathematical world congress \cite{Serfaty} and of Djalil Chafai \cite{Chafai2}.

\vspace{.2cm}

{\bf Acknowledgment.}
We thank the anonymous reviewers and  our colleagues at UCA Didier Clamond, Fran\c{c}ois Delarue and Gilles Scarella, for useful suggestions.


\end{document}